\documentclass[12pt]{article}
\usepackage{amssymb}
\usepackage{amsmath}

\newtheorem{theorem}{Theorem}
\newtheorem{corollary}{Corollary}
\newtheorem{lemma}{Lemma}

\begin{document}

\addvspace{2in}

\title{Global Existence for the Vlasov-Poisson System with Steady Spatial Asymptotics}
\author{Stephen Pankavich \\
Department of Mathematical Sciences \\
Carnegie Mellon University \\
Pittsburgh, PA 15213 \\
sdp@andrew.cmu.edu}
\date{\today}
\maketitle

\begin{center}
\textbf{MOS Classification :} 35L60, 35Q99, 82C21, 82C22, 82D10. 
\end{center}

\vspace{0.3in}

\begin{abstract}
A collisionless plasma is modelled by the Vlasov-Poisson system in
three space dimensions.  A fixed background of positive charge,
dependant upon only velocity, is assumed.  The situation in which
mobile negative ions balance the positive charge as $\vert x \vert
\rightarrow \infty$ is considered.  Thus, the total positive
charge and the total negative charge are both infinite. Smooth
solutions with appropriate asymptotic behavior for large $\vert x
\vert$, which were previously shown to exist locally in time, are
continued globally. This is done by showing that the charge
density decays at least as fast as $\vert x \vert^{-6}$.  This
paper also establishes decay estimates for the electrostatic field
and its derivatives.
\end{abstract}

\begin{center}
INTRODUCTION \\
\end{center}

Let $F : \mathbb{R}^3 \rightarrow [0,\infty)$, $f_0 : \mathbb{R}^3
\times \mathbb{R}^3 \rightarrow [0,\infty)$, and $A : [0,\infty)
\times \mathbb{R}^3 \rightarrow \mathbb{R}^3$ be given. We seek a
solution, $f : [0,\infty) \times \mathbb{R}^3 \times \mathbb{R}^3
\rightarrow [0,\infty)$ satisfying

\begin{equation}
\label{one}
\left. \begin{array}{ccc}
& & \partial_t f + v \cdot \nabla_x f - (E + A) \cdot \nabla_v f = 0,\\
\\
& & \rho(t,x) = \int (F(v) - f(t,x,v)) dv,\\
\\
& & E(t,x) = \int \rho(t,y) \frac{x - y}{\vert x - y \vert^3} dy,\\
\\
& & f(0,x,v) = f_0(x,v). \end{array} \right \}
\end{equation}

Here $F$ describes a number density of positive ions which form a
fixed background, and $f$ denotes the density of mobile negative
ions in phase space.  Notice that if $f_0(x,v) = F(v)$ and $A =
0$, then $f(t,x,v) = F(v)$ is a steady solution.  Thus, we seek
solutions for which $f(t,x,v) \rightarrow F(v)$ as $\vert x \vert
\rightarrow \infty$. It is important to notice that (\ref{one}) is
a representative problem, and that problems concerning multiple
species of ions can be treated in a similar manner.\\

Precise conditions which ensure local existence and conditions for
continuation were given in \cite{VPSSA}. Also, \emph{a priori}
bounds on $\rho$ and a quantity related to the energy were obtained
in \cite{SSAVP}. Therefore, we will work towards establishing global
existence with these bounds and using similar assumptions.  Other
work on the infinite mass case has been done in \cite{J} and
\cite{C}. \\

The case when $F(v) = 0$ and $f \rightarrow 0$ as $\vert x \vert
\rightarrow \infty$ has been studied extensively.  Smooth solutions
were shown to exist globally in time in \cite{P} and independently
in \cite{LP}. Important results prior to global existence appear in
\cite{B}, \cite{GS}, \cite{Horst}, \cite{Horst82}, and \cite{K}.
Also, the method used by \cite{P} is refined in \cite{Horst93} and
\cite{Schaeffer}. Global existence for the Vlasov-Poisson system in
two dimensions was established in \cite{OU} and \cite{W}. A complete
discussion of the literature concerning Vlasov-Poisson may be found
in \cite{Glassey}. We also mention \cite{BR} and \cite{RR} since the
problem treated in these papers is periodic in space, and thus the
solution does not decay for large $\vert x \vert$.\\

\section{Section 1}
Let $p \in (3,4)$ and denote
$$R(x) = R(\vert x \vert) = (1 + \vert x \vert^2)^\frac{1}{2}.$$
We will use the following notation :

$$\Vert g \Vert_{L^\infty(\mathbb{R}^n)} = \sup_{z \in
\mathbb{R}^n} \vert g(z) \vert $$ and

$$ \Vert \rho \Vert_p := \Vert \rho(x) R^p(x)
\Vert_{L^\infty(\mathbb{R}^3)},$$ but never use $L^p$. We will
write, for example, $\Vert \rho(t) \Vert_p$ for the $\Vert
\cdot \Vert_p$ norm of $x \mapsto \rho(t,x)$.\\

Following \cite{SSAVP} and \cite{VPSSA}, we assume the following
conditions hold for some $C > 0$ and all $t \geq 0$, $x \in
\mathbb{R}^3$, and $v \in \mathbb{R}^3$, unless otherwise stated :

\begin{enumerate}
\renewcommand{\labelenumi}{(\Roman{enumi})}
\item $F(v) = F_R(\vert v \vert)$ is nonnegative and $C^2$ with

\begin{equation}
F^{''}_R(0) < 0,
\end{equation}

and there is $W \in (0,\infty)$ such that

\begin{equation}
\label{W}
\left . \begin{array}{ccc}
F_R^\prime(u) < 0 \ \ \ \ \mathrm{for} \ u \in (0,W)\\
\\
F_R(u) = 0 \ \ \ \ \mathrm{for} \ u \geq W.
\end{array} \right \}
\end{equation}

\item $f_0$ is $C^1$ and nonnegative.

\item $A$ is $C^1$ with

\begin{equation}
\vert A(t,x) \vert \leq C^{(0)} R^{-2}(x),
\end{equation}

\begin{equation}
\vert \partial_{x_i} A(t,x) \vert \leq C R^{-3}(x),
\end{equation}

and

\begin{equation}
\nabla_x \cdot A(t,x) = 0.
\end{equation}

Finally, we assume there is a continuous function $a : [0,T]
\rightarrow \mathbb{R}$ such that
$$\left \vert A(t, x) - a(t) \frac{x}{\vert x \vert^3} \right \vert \leq C R^{1-p}(x).$$

\item $F - f_0$ has compact support in $v$, and there is $N > 0$
such that for $\vert x \vert > N$, we have
$$ \vert F(v) - f_0(x,v) \vert \leq CR^{-6}(x).$$

\end{enumerate}

\noindent Then, we have global existence:

\begin{theorem}
Assuming conditions $(I)$-$(IV)$ hold, there exists $f \in
\mathcal{C}^1([0,\infty) \times \mathbb{R}^3 \times \mathbb{R}^3)$
that satisfies (\ref{one}) with $\Vert \int (F - f)(t) \ dv
\Vert_p$ bounded on $t \in [0,T]$, for every $T > 0$. Moreover,
$f$ is unique.
\end{theorem}

\vspace{0.3in}

In addition, due to the previously known result of \cite{Pankavich}, stated here as Theorem $2$, we are able to conclude further decay of the charge density in Corollary $1$ :

\begin{theorem}
Let $T > 0$ and $f$ be the \ $\mathcal{C}^1$ solution of (\ref{one}) on $[0,T] \times
\mathbb{R}^3 \times \mathbb{R}^3$.  Then, we have

$$ \Vert \rho(t) \Vert_6 \leq C_{p, t}  $$
for any $t \in [0,T]$, where $C_{p, t}$ depends upon
$$\sup_{\tau \in [0,t]} \Vert \rho(\tau) \Vert_p.$$
\end{theorem}

\begin{corollary}
Let $T > 0$ and $f$ be the \ $\mathcal{C}^1$ solution of (\ref{one}) on $[0,T] \times
\mathbb{R}^3 \times \mathbb{R}^3$.  Then, we have
$$\left \Vert \int (F - f)(t) \ dv \right \Vert_6 \leq C$$ for any $t \in [0,T]$.
\end{corollary}

\vspace{0.3in}

To prove Theorem $1$, we will use Theorems $1$ and $3$ of
\cite{VPSSA}, which guarantee local existence and continuation of
the local solution so long as $\Vert \int (F - f)(t) \Vert_p$ is
bounded for some $p > 3$. Therefore, the following lemma will be
all that is needed to complete the proof of the theorem.

\begin{lemma}
Assume conditions $(I)$-$(IV)$ hold. Let $f$ be a $\mathcal{C}^1$
solution of (\ref{one}) on $[0,T] \times \mathbb{R}^3 \times
\mathbb{R}^3$. Then,

$$ \left \Vert \int (F - f)(t) \ dv \right \Vert_p \leq C $$
for all $t \in [0,T]$, where $C$ is determined by $F$, $A$, $f_0$,
and $T$.
\end{lemma}

\vspace{0.3in}

Define
\begin{equation}
\label{Q}
Q_f(t) := \sup \{\vert v \vert : \exists x \in \mathbb{R}^3,
\tau \in [0,t] \ \mathrm{such \ that} \ f(\tau, x, v) \neq 0 \}
\end{equation}
and
$$ Q_g(t) : = \max \{ W, Q_f(t) \}.$$

In Section 2, we will bound $E$, $\nabla E$, and $\nabla_{x,v} f$.
Also, we will estimate the energy, and obtain bounds on $Q_f(t)$
and thus $Q_g(t)$.  To better reveal the line of thought, the
proofs of Lemmas $2$ through $4$ are deferred to Section $3$.

We will denote by ``$C$'' a generic constant which changes from
line to line and may depend upon $f_0$, $A$, $F$, or $T$, but not
on $t$, $x$, or $v$. When it is necessary to refer to a specific
constant, we will use numeric superscripts to distinguish them.
For example, $C^{(0)}$, as in $(III)$, will always refer to the
same numerical constant.\\

To estimate $E$ and $\nabla E$, we will use
\begin{lemma}
For any $q > 0$ and $b \in [0, \frac{5}{18})$ with $b \leq \frac{2}{q}$, we have

$$ \vert E(t,x) \vert \leq C \left ( \Vert \rho(t) \Vert_q R^{-q}(x) \right )^b$$
for any $t \in [0,T]$, $\vert x \vert \geq 1$, where $C$ may depend upon $\Vert \rho(t) \Vert_\infty$.
\end{lemma}

\noindent and

\begin{lemma}
For any $q > 0$ and $a \in [0,1)$ with $a \leq \frac{3}{q}$, we have

$$ \vert \nabla E(t,x) \vert \leq C \left (\Vert \rho(t) \Vert_q R^{-q}(x) \right )^a$$
for any $t \in [0,T]$, $\vert x \vert \geq 1$, where $C$ may depend upon $\Vert \rho(t) \Vert_\infty$ and $\Vert \nabla \rho(t) \Vert_\infty$.
\end{lemma}

\vspace{0.3in}

Then, we will use the following lemma to bound the $p$-norm of
$\rho$ and obtain decay of both $E$ and $\nabla E$ :

\begin{lemma}
For any $q \in [0,\frac{54}{13})$, $\Vert \rho(t) \Vert_q$ is bounded for $t \in [0,T]$, where $C$ may depend upon $\Vert \rho(t) \Vert_\infty$, $\Vert \nabla \rho(t) \Vert_\infty$, and $Q_f(T)$.
\end{lemma}

Thus, once we show $\Vert \rho(t) \Vert_\infty$ and $\Vert \nabla \rho(t) \Vert_\infty$ are bounded for $t \in [0,T]$, and $Q_f(T)$ is finite, we may use Lemma $4$ to
bound $\Vert \rho(t) \Vert_p$ since $p \in (3,4)$.  Then, we can prove Lemma $1$, and thus
Theorem $1$.

\vspace{0.3in}

\section{Section 2}
\subsection{Characteristics}
Define the characteristics, $X(s,t,x,v)$ and $V(s,t,x,v)$, by
\begin{equation}
\label{char} \left. \begin{array}{ccc} & &
\frac{\partial X}{\partial s} (s,t,x,v) = V(s,t,x,v) \\
\\
& & \frac{\partial V}{\partial s} (s,t,x,v) = - \left( E(s,X(s,t,x,v)) + A(s, X(s,t,x,v)) \right) \\
\\
& & X(t,t,x,v) = x \\
\\
& & V(t,t,x,v) = v. \end{array} \right \}
\end{equation}

\noindent Then, we have
$$ \frac{\partial}{\partial s} f(s,X(s,t,x,v),V(s,t,x,v)) =
\partial_t f + V \cdot \nabla_x f - ( E + A ) \cdot \nabla_v f = 0.$$

\noindent Therefore, $f$ is constant along characteristics, and
 \begin{equation}
\label{fchar} f(t,x,v) = f(0, X(0,t,x,v), V(0,t,x,v)) =
f_0(X(0,t,x,v),V(0,t,x,v)).
\end{equation}
Thus, we find by $(II)$ that $f$ is nonnegative.\\

\noindent Define $g : [0, \infty) \times \mathbb{R}^3 \times
\mathbb{R}^3 \rightarrow \mathbb{R}$ by $$ g(t,x,v) := F(v) -
f(t,x,v) .$$ Then, we see that $ \sup_{x,v} \vert g \vert \leq \Vert F \Vert_{L^\infty} + \Vert f_0
\Vert_{L^\infty} < \infty$, and
\begin{equation}
\label{gchar}
\begin{array}{rcl}
\frac{\partial}{\partial s} g(s,X(s,t,x,v),V(s,t,x,v)) & = &
\partial_t g + V \cdot \nabla_x g - ( E + A ) \cdot \nabla_v g \\
\\
& = & - \nabla_v F(V(s)) \cdot (E + A)(s,X(s)).
\end{array}
\end{equation}
Finally, for any $s \in [0,T]$ with $f(s,X(s),V(s)) \neq 0$,
we have for $t \in [0,T]$ and $x,v \in \mathbb{R}^3$,
\begin{eqnarray*}
\vert X(s,t,x,v) - x \vert & = & \left \vert \int_s^t \dot{X}(\tau,t,x,v) \ d\tau \right \vert \\
& \leq & \int_s^t \vert V(\tau,t,x,v) \vert \ d\tau \\
& \leq & TQ_g(T).
\end{eqnarray*}

\noindent So, assuming we can bound $Q_g(T)$, we have for $\vert x
\vert \geq 2TQ_g(T)$, $v \in \mathbb{R}^3$, $t \in [0,T]$, and $s
\in [0,t]$,

\begin{equation}
\label{xchar1} \vert X(s,t,x,v) \vert \geq \vert x \vert - T Q_g(T)
\geq \frac{1}{2} \vert x \vert
\end{equation}
and
\begin{equation}
\label{xchar2} \vert X(s,t,x,v) \vert \leq \vert x \vert + T Q_g(T)
\leq \frac{3}{2} \vert x \vert.
\end{equation}

Unless it is necessary, we will omit writing the dependence of
$X(s)$ and $V(s)$ on
$t$, $x$, and $v$ for the remainder of the paper.\\

\vspace{0.15in}

\subsection{Bounds on the Field}

Much of the work that follows will rely on energy estimates found
in \cite{SSAVP}.  In particular, we combine Lemma $3$ and Theorem
$4$ from \cite{SSAVP} to obtain an a priori bound on $\rho$.\\

First, define $\sigma : [0,\infty) \rightarrow \mathbb{R}$ by
$$ \sigma(h) = - \int_0^{\min \{h, F(0) \}} \left ( F^{-1}(\tilde{h})
\right )^2 \ d\tilde{h} $$ and $S : [0,\infty) \times [0,\infty)
\rightarrow [0,\infty) $ by
$$ S(h, \eta) = (h - F(\eta))\eta^2 + \sigma(h) -
\sigma(F(\eta)).$$

Then, from \cite{SSAVP}, we have the following lemma:

\begin{lemma}
Let $$k(t,x) = \int S \left (f(t,x,v), \vert v \vert \right ) \ dv
.$$ Assuming conditions $(I) - (IV)$ hold, we have
\begin{equation}
\label{conslaw1}
\int \vert F(v) - f(t,x,v) \vert \ dv \leq
C(k(t,x)^\frac{3}{5} + k(t,x)^\frac{1}{2})
\end{equation}
and
\begin{equation}
\label{conslaw2} \int k(t,x) \ dx \leq C.
\end{equation}
\end{lemma}

\vspace{0.3in}

It follows directly from (\ref{conslaw1}) that
\begin{equation}
\label{conslaw3}
\vert \rho(t,x) \vert \leq C (k(t,x)^\frac{3}{5}
+ k(t,x)^\frac{1}{2})
\end{equation}

\noindent Then, using Lemma $5$ , we may bound the field.\\

Assume $\Vert \rho(t) \Vert_{L^\infty(\mathbb{R}^6)} \leq C$ for
all $t \in [0,T]$. Then, for $t \geq 0$, $x \in \mathbb{R}^3$, and
any $R > 0$
\begin{equation}
\label{Emethod}
\left.
\begin{array}{rcl}
\vert E(t,x) \vert & \leq & \int \frac{\vert \rho(t,y)
\vert}{\vert x - y \vert^2} \ dy \\
\\
& \leq & \int_{\vert x - y \vert < R} \Vert \rho(t)
\Vert_{L^\infty} \vert x - y \vert^{-2} \ dy + \int_{\vert x - y
\vert > R} \vert \rho(t,y) \vert \ \vert x -y \vert^{-2} \ dy \\
\\
& \leq & 4 \pi \Vert \rho(t) \Vert_{L^\infty} \int_0^R \ dr +
C \int_{\vert x - y \vert > R} (k^\frac{1}{2} (t,y) +
k^\frac{3}{5}(t,y)) \vert x - y \vert^{-2} \ dy \\
\\
& \leq & 4 \pi \Vert \rho(t) \Vert_{L^\infty} R + C (\int k(t,y) \
dy)^\frac{1}{2} (\int_{\vert x - y \vert > R} \vert x - y
\vert^{-4})^\frac{1}{2} \\
\\
& & \ \ \ + C (\int k(t,y) \ dy)^\frac{3}{5} (\int_{\vert x - y
\vert > R} \vert x - y \vert^{-5})^\frac{2}{5} \\
\\
& \leq & 4 \pi \Vert \rho(t) \Vert_{L^\infty} R + C (\int_R^\infty
r^{-2} \ dr)^\frac{1}{2} + C(\int_R^\infty r^{-3} \
dr)^\frac{2}{5} \\
\\
& \leq & C(\Vert \rho(t) \Vert_{L^\infty} R + R^{-\frac{1}{2}} +
R^{-\frac{4}{5}}).
\end{array} \right \}
\end{equation}

Obviously, we may choose $R = 1$ and deduce that $\vert E(t,x)
\vert \leq C$.  Suppose that the $v$-support of $g$ is bounded for
bounded times. Then, we find

$$ \vert \rho(t,x) \vert \leq \int_{\vert v \vert \leq Q_g(t)} \vert
g(t,x,v) \vert \ dv \leq \left ( \Vert f_0 \Vert_{L^\infty}  + \Vert F \Vert_{L^\infty} \right ) Q_g^3(t) \leq C
Q_g^3(t).$$

Thus, we know for every $t \in [0,T]$,

$$ \Vert \rho(t) \Vert_{L^\infty} \leq C Q_g^3(t). $$

Now, choose $R = Q_g^{-\frac{5}{3}} (t)$.  If $R \geq 1$, then
$Q_g(t)$ is bounded.  Otherwise, $R \leq 1$ and thus $
R^{-\frac{1}{2}} \leq R^{-\frac{4}{5}}$.

Then, we have
\begin{equation}
\label{E2method}
\left.
\begin{array}{rcl}
\vert E(t,x) \vert & \leq & C(\Vert \rho(t) \Vert_{L^\infty} R +
R^{-\frac{4}{5}}) \\
\\
& \leq & C \left ( Q_g^3(t) Q_g^{-\frac{5}{3}}(t) +
(Q_g^{-\frac{5}{3}}(t))^{-\frac{4}{5}} \right ) \\
\\
& \leq & C Q_g^\frac{4}{3} (t) \\
\\
& =: & C^{(1)} Q_g^\frac{4}{3} (t)
\end{array}
\right \}
\end{equation}

\subsection{Bounds on Derivatives of Density and Field}
We proceed as in Section $4.2.5$ of \cite{Glassey} with one modification :  we will
not assume that $\rho(t) \in L^1(\mathbb{R}^3)$ for all $t \in [0,T]$.  Instead, we find
for any $R > 0$, (letting $r = \vert x - y \vert$)

\begin{eqnarray*}
\left \vert \frac{1}{4\pi} \int_{\vert y - x \vert \geq R}
\rho(t,y) \left ( \frac{3(y_k - x_k)^2}{r^5} - \frac{1}{r^3}
\right ) \ dy \right \vert & \leq & \frac{1}{\pi} \int_{\vert y -
x \vert \geq R} \vert \rho(t,y) \vert \ r^{-3} \ dy \\
& \leq & C \int_{\vert y - x \vert \geq R} \left (
k^\frac{1}{2}(t,y) + k^\frac{3}{5}(t,y) \right ) r^{-3} \ dy \\
& \leq & C \left ( (\int k(t,y) \ dy)^\frac{1}{2} (\int_{\vert y -
x \vert \geq R} r^{-6} \ dy)^\frac{1}{2} \right. \\
& \ & \ \left. + (\int k(t,y) \ dy)^\frac{3}{5} (\int_{\vert y - x \vert \geq R} r^{-\frac{15}{2}} \ dy)^\frac{2}{5} \right ) \\
& \leq & C \left ( R^{-\frac{3}{2}} + R^{-\frac{9}{5}} \right ).
\end{eqnarray*}

Thus, introducing this estimate into the result of \cite{Glassey},
we have for any $0 < d \leq R$
\begin{equation}
\label{derivfield} \vert \nabla_x E(t,x) \vert \leq C(1 +
\ln(R/d))\Vert \rho(t) \Vert_{L^\infty} + C d \Vert \nabla_x
\rho(t) \Vert_{L^\infty} + C ( R^{-\frac{3}{2}} +
R^{-\frac{9}{5}}).
\end{equation}

\noindent Then, we may follow the argument in Section $4.2.6$ of \cite{Glassey} and
find a priori bounds on $\Vert \nabla_x E(t) \Vert_{L^\infty}$
and $\Vert \nabla_{x,v} f(t) \Vert_{L^\infty}$ as long as $\Vert
\rho(t) \Vert_{L^\infty}$ is bounded.  This work will be included in Appendix A.

\subsection{Energy Bound}

First, notice from (\ref{W}) that we may write $W > 0$ as $$W = \inf \{
\eta > 0 : F(\eta) = 0 \}.$$
Recall the definitions of $\sigma$, $S$, and $k$, as well as,
(\ref{conslaw2}), which we may write as
$$ \int \int S(f(t,x,v), \vert v \vert) \ dv \ dx \leq C $$ for any
$t \in [0,T]$.

Let $\eta \geq 2W$. We find
\begin{eqnarray*}
S(h,\eta) & = & (h - F(\eta))\eta^2 + \sigma(h) - \sigma(F(\eta)) \\
& = & h \eta^2 + \sigma(h) - \sigma(0) \\
& = & h \eta^2 + \sigma(h) \\
& = & h \eta^2 - \int_0^{\min\{h, F(0)\}} (F^{-1}(\tilde{h}))^2 d\tilde{h} \\
& \geq & h \eta^2 - \min\{h, F(0)\} (F^{-1}(0))^2 \\
& \geq & h \eta^2 - hW^2 \\
& = & h(\eta^2 - W^2).
\end{eqnarray*}

Then, since $\eta \geq 2W$, we have $\frac{1}{4} \eta^2 \geq W^2$
and
\begin{eqnarray*}
h (\eta^2 - W^2) & = & \frac{1}{2}h \eta^2 + h(\frac{1}{2}\eta^2 -
W^2) \\
& \geq & \frac{1}{2} h\eta^2.
\end{eqnarray*}

Thus, for any $h \geq 0$ and $\eta \geq 2W$,
$$ S(h,\eta) \geq \frac{1}{2} h \eta^2 $$

and, finally, for $P \geq 2W$

\begin{equation}
\label{energy} \int \int_{\vert v \vert > P} \vert v \vert^2
f(t,x,v) \ dv \ dx \leq 2 \int \int S(f(t,x,v), \vert v \vert) \
dv \ dx \leq C.
\end{equation}

\subsection{The Good, the Bad, and the Ugly Revisited}
The method we will employ is very similar to that used in Section $4.4$ of
\cite{Glassey}.  The differences are due mainly to the lack of
positivity in the charge density, changes in the
conserved energy, and the contribution of the applied field.  We will
assume throughout this section that $Q_g(t)$ is not already bounded and
(\ref{E2method}) applies.

Define $$Q(t) := \left ( \max \{(2W)^\frac{4}{3}, C^{(0)}
\}\right )^\frac{15}{13} + Q_f(t).$$
Then, define
$C^{(2)} := (C^{(0)})^{-\frac{7}{13}}
+ C^{(1)}$ so that we use (\ref{E2method}) and $(III)$  to find
\begin{eqnarray*}
\vert E(\tau,x) + A(\tau,x) \vert & \leq & C^{(1)}Q^\frac{4}{3}(t)
+ C^{(0)} R^{-2}(x) \\
& \leq & C^{(1)}Q^\frac{4}{3}(t) + \left ( C^{(0)} \right )^{-\frac{7}{13}}
\left ( C^{(0)} \right )^\frac{20}{13} \\
& \leq & C^{(1)}Q^\frac{4}{3}(t) + \left ( C^{(0)} \right
)^{-\frac{7}{13}} Q^\frac{4}{3}(t) \\
& = & C^{(2)} Q^\frac{4}{3}(t)
\end{eqnarray*}
for any $\tau \in [0,t]$.\\

Now, let $(\hat{X}(t), \hat{V}(t))$ be any fixed characteristic :
$$ \frac{d}{dt} \hat{X} = \hat{V}, \frac{d}{dt} \hat{V} =
E(t,\hat{X}) $$ for which
$$ f(t,\hat{X}(t), \hat{V}(t)) \neq 0.$$

For any $0 \leq \Delta \leq t$, we have
\begin{eqnarray}
\int_{t - \Delta}^t \vert E(s,\hat{X}(s)) \vert \ ds & \leq & C
\int_{t - \Delta}^t \int \int \frac{ \vert g(s,y,w) \vert}{\vert y - \hat{X}(s) \vert^2}
\ dw \ dy \ ds \\
\label{intE} & = & C \int_{t - \Delta}^t \int \int \frac{\vert
g(s,X(s,t,x,v),V(s,t,x,v)) \vert}{\vert X(s,t,x,v) - \hat{X}(s)
\vert^2} \ dv \ dx \ ds.
\end{eqnarray}

Let $P = Q^\frac{13}{20}(t)$, $R > 0$, and $\Delta =
\frac{P}{4C^{(2)}Q^\frac{4}{3}(t)}$. From the definition of $Q$,
notice that $P \geq 2W$.\\

Let us partition the integral into $I_G$, $I_B$, and $I_U$, where
$I_A$ is the integral in (\ref{intE}) over the set $A$, and the three sets are defined as :
\begin{eqnarray*}
G & := & \{(s,x,v) : t - \Delta < s < t \ \mathrm{and} \ (\vert v \vert < P
\ \mathrm{or} \ \vert v - \hat{V}(t) \vert < P) \}, \\
\\
B & := & \{(s,x,v) : t - \Delta < s < t \ \mathrm{and} \  \vert v \vert > P \ \mathrm{and} \
\vert v - \hat{V}(t) \vert > P \\
& & \ \mathrm{and} \ \vert X(s,t,x,v) - \hat{X}(s) \vert < R \}, \\
\\
U & := & \{(s,x,v) : t - \Delta < s < t \ \mathrm{and} \ \vert v \vert > P \ \mathrm{and} \
\vert v - \hat{V}(t) \vert > P \\
& & \ \mathrm{and} \ \vert X(s,t,x,v) - \hat{X}(s) \vert > R \}.
\end{eqnarray*}

We will use the invertibility of the characteristics as described
in \cite{Glassey}, so that when we set
$$ y = X(s,t,x,v) $$
$$ w = V(s,t,x,v) $$
we can invert using
$$ x = X(t,s,y,w) $$
$$ v = V(t,s,y,w). $$
In particular, notice $w = V(s,t,X(t,s,y,w),V(t,s,y,w))$ and $\frac{\partial(y,w)}{\partial(x,v)} = 1$ .\\

To handle the integral over $G$, we must first deal with some
preliminary inequalities :

\begin{enumerate}
\item First notice that $\vert V(s,t,x,v) - v \vert \leq \int_s^t
\left \vert E(\tau,X(\tau)) + A(\tau,X(\tau)) \right \vert \
d\tau$.

Thus, for $s \in [t - \Delta, t]$, we have
$$ \vert V(s,t,x,v) - v \vert \leq \Delta C^{(2)} Q^\frac{4}{3}(t) = \frac{1}{4}P.$$

\item For $\vert v \vert < P$,
$$\vert V(s,t,x,v) \vert \leq
\vert v \vert + \frac{1}{4}P < 2P.$$

\item For $\vert v - \hat{V}(t) \vert < P$,
$$\vert V(s,t,x,v) -
\hat{V}(s) \vert \leq \vert v - \hat{V}(t) \vert + \vert
\hat{V}(t) - \hat{V}(s) \vert + \vert V(s,t,x,v) - v \vert \leq P
+ \frac{1}{4}P + \frac{1}{4}P < 2P.$$

\item For $\vert v \vert > P$ and $\tau \in [t - \Delta, t]$,
\begin{eqnarray*}
\vert V(\tau,t,x,v) \vert & \geq & \vert v \vert - \frac{1}{4}P \\
& > & \frac{3}{4} P \\
& \geq & \frac{3}{4} \cdot 2W \\
& > & W
\end{eqnarray*}

\end{enumerate}

Now, let $$\chi_G(s,x,v) := \left \{ \begin{array}{lc} 1 & (s,x,v)
\in G \\ 0 & $ else. $ \end{array} \right. $$ Then, we have
\begin{eqnarray*}
I_G & = & \mathop{\int \int \int}_G \frac{\vert
g(s,X(s,t,x,v),V(s,t,x,v)) \vert}{\vert X(s,t,x,v) - \hat{X}(s)
\vert^2}
\ dv \ dx \ ds \\
& = & \int_{t - \Delta}^t \int \int \frac{\chi_G(s,x,v) \vert
g(s,X(s,t,x,v),V(s,t,x,v)) \vert}{\vert X(s,t,x,v) - \hat{X}(s) \vert^2} \ dv \ dx \ ds \\
& = & \int_{t - \Delta}^t \int \int
\frac{\chi_G(s,X(t,s,y,w),V(t,s,y,w)) \vert g(s,y,w) \vert}{\vert
y - \hat{X}(s) \vert^2} \ dw \ dy \ ds.
\end{eqnarray*}
If $\chi_G \neq 0$, then $\vert V(t,s,y,w) \vert < P $ or $\vert
V(t,s,y,w) - \hat{V}(t) \vert < P$.  Then, by Preliminary Inequalities $2$ and $3$, we
have either $\vert w \vert < 2P$ or $\vert w - \hat{V}(s) \vert <
2P$. Set $$ \tilde{\rho}(s,y) := \int \vert g(s,y,w) \vert
\chi_G(s,X(t,s,y,w), V(t,s,y,w)) \ dw.$$ Then, $\Vert
\tilde{\rho}(s) \Vert_{L^\infty} \leq C P^3$. Also, using
(\ref{conslaw1}), we know $$\tilde{\rho}(s,y) \leq \int \vert
g(s,y,w) \vert \ dw \leq C(k^\frac{3}{5}(s,y) +
k^\frac{1}{2}(s,y)).$$  Thus, we employ the method of (\ref{Emethod}) and (\ref{E2method}) to find

$$ \int \frac{\tilde{\rho}(s,y)}{\vert y - \hat{X}(s) \vert^2} \ dy
\leq C P^\frac{4}{3}$$ and, finally, we have
\begin{equation}
\label{IG} I_G = \int_{t - \Delta}^t \int
\frac{\tilde{\rho}(s,y)}{\vert y - \hat{X}(s) \vert^2} \ dy \ ds
\leq C \Delta P^\frac{4}{3}.
\end{equation}

Estimating $I_B$, we have
\begin{eqnarray*}
I_B & = & \mathop{\int \int \int}_B \frac{\vert g(s,X(s,t,x,v),V(s,t,x,v)) \vert}
{\vert X(s,t,x,v) - \hat{X}(s) \vert^2} \ dv \ dx \ ds \\
& \leq & \int_{t - \Delta}^t \int_{\vert y - \hat{X}(s) \vert <
R} \int \frac{\vert g(s,y,w) \vert}{\vert y - \hat{X}(s)
\vert^2} \ dw \ dy \ ds \\
& \leq & \left ( \Vert F \Vert_{L^\infty} + \Vert  f_0 \Vert_{L^\infty} \right ) Q^3(t) \int_{t - \Delta}^t
\int_{\vert y - \hat{X}(s) \vert < R} \vert y - \hat{X}(s)
\vert^{-2} \ dw \ ds \\
& \leq & C Q^3(t) \int_{t - \Delta}^t \int_0^R \ dr \ ds,
\end{eqnarray*}
and thus

\begin{equation}
\label{IB}
I_B \leq C \Delta Q^3(t) R.
\end{equation}

Finally, to estimate $I_U$, we use Section $3$ (specifically, line $(15)$) of \cite{Schaeffer} to find
\begin{equation}
\label{uglyeq}
\int_{t - \Delta}^t \vert X(s,t,x,v) - \hat{X}(s)
\vert^{-2} \ \chi_U(s,x,v) \ ds \leq \frac{C}{RP},
\end{equation}
for $(x,v)$ as in $U$.
The proof of this result is quite long, so we shall include a sketch rather than the entire proof.  First, let
$$ Z(s) = X(s,t,x,v) - \hat{X}(s).$$
Then, choose $s_0 \in [t - \Delta, t]$ such that $$ \vert Z(s_0) \vert \leq \vert Z(s) \vert$$
for all $s \in [t - \Delta, t]$.  It is shown in \cite{Schaeffer} that
\begin{equation}
\label{Zineq}
\vert Z(s) \vert \geq \frac{1}{4} P \vert s - s_0 \vert.
\end{equation}
Now, define $$ \Sigma(r) := \left \{ \begin{array}{cc} \frac{1}{R^2} & 0 \leq r \leq R^2 \\ \ & \ \\ \frac{1}{r} & r \geq R^2. \end{array} \right.$$
Notice that $\Sigma$ is non-negative, non-increasing and $$ \vert Z(s) \vert^{-2} \chi_U(s,x,v) \leq \Sigma(\vert Z(s) \vert^2).$$
Using these properties of $\Sigma$ with (\ref{Zineq}), we find
\begin{eqnarray*}
\int_{t - \Delta}^t \vert Z(s) \vert^{-2} \chi_U(s,x,v) \ ds & \leq & \int_{t - \Delta}^t \Sigma(\vert Z(s) \vert^2) \ ds \\
& \leq & \int_{t - \Delta}^t \Sigma \left ( (\frac{1}{4}P \vert s - s_0 \vert)^2 \right ) \ ds \\
& \leq & \int \Sigma \left ( \frac{P^2}{16} \tau^2 \right ) \ d\tau \\
& = & \frac{16}{PR}.
\end{eqnarray*}
This shows (\ref{uglyeq}).

Now, using (\ref{uglyeq}) with (\ref{W}), (\ref{fchar}), (\ref{energy}), and preliminary inequality $4$, we have
\begin{eqnarray*}
I_U & = & \mathop{\int \int \int}_U \frac{\vert
g(s,X(s,t,x,v),V(s,t,x,v)) \vert}{\vert X(s,t,x,v) - \hat{X}(s) \vert^2} \ dv \ dx \ ds \\
& = & \mathop{\int \int \int}_U \frac{\vert F(V(s,t,x,v)) - f(s, X(s,t,x,v), V(s,t,x,v)) \vert}
{\vert X(s,t,x,v) - \hat{X}(s) \vert^2} \ dv \ dx \ ds \\
& = & \mathop{\int \int \int}_U \frac{\vert - f(t,x,v) \vert}
{\vert X(s,t,x,v) - \hat{X}(s) \vert^2} \ dv \ dx \ ds \\
& = & \int_{t - \Delta}^t \int_{\vert v \vert > P \ \cap \ \vert v
- \hat{V}(t) \vert > P} \int_{\vert X(s,t,x,v) - \hat{X}(s)
\vert > R}
\frac{f(t,x,v)} {\vert X(s,t,x,v) - \hat{X}(s) \vert^2} \ dx \ dv \ ds \\
& = & \int_{\vert v \vert > P} \int f(t,x,v) \
\left ( \int_{t - \Delta}^t \vert X(s,t,x,v) - \hat{X}(s) \vert^{-2} \  \chi_U(s,x,v) \ ds \right ) \ dx \ dv  \\
& \leq & \frac{C}{RP} \int \int_{\vert v \vert > P} f(t,x,v) \ dv \ dx \\
& \leq & \frac{C}{RP^3} \int \int_{\vert v \vert > P} \vert v
\vert^2 \ f(t,x,v) \ dv \ dx
\end{eqnarray*}
and so

\begin{equation}
\label{IU}
I_U \leq \frac{C}{RP^3}.
\end{equation}

Finally, collecting the estimates (\ref{IG}), (\ref{IB}), and
(\ref{IU}), we find,

$$ \frac{1}{\Delta} \int_{t - \Delta}^t \vert E(s,\hat{X}(s)) \vert \
ds \leq C \left( P^\frac{4}{3} + R Q^3(t)  + \frac{1}{\Delta RP^3}
\right ).
$$ We take $R = Q^{-\frac{32}{15}}(t)$
and then
$$ \frac{1}{\Delta} \int_{t - \Delta}^t \vert E(s,\hat{X}(s)) \vert \
ds \leq C Q^\frac{13}{15}(t). $$

\noindent Using $(III)$, we have
\begin{equation}
\label{GBU}
\begin{array}{rcl}
\frac{1}{\Delta} \int_{t - \Delta}^t \left ( \vert E(s,\hat{X}(s))
+ A(s, \hat{X}(s)) \vert \right ) \ ds & \leq & C
Q^\frac{13}{15}(t) + \frac{1}{\Delta} \int_{t-\Delta}^t
C^{(0)} R^{-2}(\hat{X}(s)) \ ds \\
& \leq & C Q^\frac{13}{15}(t) + \frac{1}{\Delta} \int_{t-\Delta}^t C^{(0)} \ ds \\
& \leq & C Q^\frac{13}{15}(t) + \left ( (C^{(0)})^\frac{15}{13}
\right
)^\frac{13}{15} \\
& \leq & C Q^\frac{13}{15}(t) + Q^\frac{13}{15}(t) \\
& \leq & C Q^\frac{13}{15}(t).
\end{array}
\end{equation}

\noindent Finally, we use the argument in Section $4.5$ of \cite{Glassey} to bound
the velocity support, since the power of $Q(t)$ is less than one.  This work will be explored in Appendix B.
Thus, for all $t \in [0.T]$,
\begin{equation}
\label{vsupport}
Q(t) \leq C,
\end{equation}
and this implies bounds on $Q_f(t)$ and $Q_g(t)$ for all $t \in
[0,T]$. Furthermore, if $V(s,t,x,v)$ satisfies
$f_0(X(0,t,x,v),V(0,t,x,v)) \neq 0$, then
$$ \vert V(s,t,x,v) \vert \leq C$$
for any $s \in [0,T]$, including $V(t,t,x,v) = v$.\\

\noindent Notice then, the bound on the velocity support implies a
priori bounds on $\Vert \rho(t) \Vert_{L^\infty}$, $\Vert E(t)
\Vert_{L^\infty}$, $\Vert \nabla_x E(t) \Vert_{L^\infty}$, and
$\Vert \nabla_{x,v} f(t) \Vert_{L^\infty}$ for all $t \in
[0,T]$.\\

Now that $\Vert \rho(t) \Vert_\infty \leq C$ and $\Vert \nabla \rho(t) \Vert_\infty \leq C$
for all $t \in [0,T]$, and $Q_f(T)$ is finite, we apply Lemma $4$ since $p > 3$, and find
$\Vert \rho(t) \Vert_p \leq C$ for all $t \in [0,T]$, and the proof of Theorem $1$ is complete.

\section{Section 3}

To conclude the paper, this section contains the proofs of Lemmas
$2$ through $4$.\\

\noindent {\bf \underline{Proof} (Lemma $2$) :} \ Let $q, T > 0$
be given with $\Vert \rho(t) \Vert_\infty \leq C$ for all $t \in
[0,T]$.  Let $b \in [0,\frac{5}{18})$ be given with $b \leq
\frac{2}{q}$. Consider $\vert x \vert \geq 1$ and define

$$ \eta : = \left ( \Vert \rho(t) \Vert_q R^{-q}(x) \right )^b $$
and divide the field into the following pieces :

$$ \vert E(t,x) \vert \leq I + II + III .$$
where

$$ I := \int_{\vert x - y \vert < \eta} \vert \rho(t,y) \vert \ \vert x - y \vert^{-2} \ dy,$$

$$ II := \int_{\eta < \vert x - y \vert < \frac{1}{2}\vert x \vert} \vert \rho(t,y) \vert \ \vert x - y \vert^{-2} \ dy, $$

$$ III := \int_{\vert x - y \vert > \frac{1}{2}\vert x \vert} \vert \rho(t,y) \vert \ \vert x - y \vert^{-2} \ dy. $$
Then, the first estimate satisfies
\begin{equation}
I \leq  C \Vert \rho(t) \Vert_\infty \eta.
\end{equation}

The second estimate satisfies, for any $m \in [0,\frac{1}{3})$,

\begin{eqnarray*}
II & \leq & C \left ( \Vert \rho(t) \Vert_q R^{-q}(x) \right )^m
\int_{\eta < \vert x - y \vert < \frac{1}{2}\vert x \vert} \vert
\rho(t,y) \vert^{1-m} \vert x - y \vert^{-2} \ dy \\
& \leq & C \left ( \Vert \rho(t) \Vert_q R^{-q}(x) \right )^m
\int_{\eta < \vert x - y \vert < \frac{1}{2}\vert x \vert} \left
(k^\frac{3(1-m)}{5}(t,y) \vert x - y \vert^{-2} +
k^\frac{1-m}{2}(t,y) \vert x - y \vert^{-2} \right ) \ dy \\
& \leq & C \left ( \Vert \rho(t) \Vert_q R^{-q}(x) \right )^m
\left [ (\int_{\eta < \vert x - y \vert < \frac{1}{2}\vert x
\vert} \vert x - y \vert^{-\frac{4}{1+m}} \ dy )
^\frac{1+m}{2} \right. \\
& \ & \ \ \left. + \ (\int_{\eta < \vert x - y \vert < \frac{1}{2}\vert x \vert}
\vert x - y \vert^{-\frac{10}{3m+2}} \ dy)^\frac{3m+2}{5} \right ] \\
& \leq & C \left ( \Vert \rho(t) \Vert_q R^{-q}(x) \right )^m \left [ \eta^{-2 + 3(\frac{1+m}{2})} +
\eta^{-2 + 3(\frac{3m+2}{5})} \right ] \\
& \leq & C \left (\Vert \rho(t) \Vert_q R^{-q}(x) \right )^m
\left \{ \begin{array}{cc}  \eta^\frac{3m-1}{2}, & \eta \geq 1 \\
\eta^\frac{9m-4}{5}, & \eta \leq 1  \end{array} \right.
\end{eqnarray*}

Then, for $\eta \geq 1$, we choose $m = \frac{3b}{2 + 3b}$, and for $\eta \leq 1$, we choose
$m = \frac{9b}{5 + 9b}$.  To guarantee convergence of the above integrals, we must have $m < \frac{1}{3}$,
and thus, $b < \frac{5}{18}$.
Thus, we find

\begin{equation}
II \leq  C \left (\Vert \rho(t) \Vert_q R^{-q}(x) \right )^b.
\end{equation}

Finally, for $b \leq \frac{2}{q}$,
\begin{eqnarray*}
III & \leq & \Vert \rho(t) \Vert_q^b (\frac{1}{2} \vert x
\vert)^{-qb} \int_{\vert x - y \vert >
\frac{1}{2} \vert x \vert} \vert \rho(t,y) \vert^{1-b} \ \vert x - y \vert^{bq - 2} \ R^{-bq}(y) \ dy \\
& \leq & C \left (\Vert \rho(t) \Vert_q R^{-q}(x) \right )^b \int_{\vert x - y \vert >
\frac{1}{2} \vert x \vert} \vert \rho(t,y) \vert^{1-b} \ (\frac{1}{4} R(y))^{bq-2} R^{-bq}(y) \ dy \\
& \leq & C \left (\Vert \rho(t) \Vert_q R^{-q}(x) \right )^b \int_{\vert x - y \vert >
\frac{1}{2} \vert x \vert} (k^\frac{1-b}{2} (t,y) + k^{\frac{3}{5}(1-b)}(t,y)) R^{-2}(y) \ dy \\
& \leq & C \left (\Vert \rho(t) \Vert_q R^{-q}(x) \right )^b \left [(\int R^{-\frac{4}{1+b}}(y) \ dy)^\frac{1+b}{2}
+ (\int R^{-\frac{10}{3b+2}}(y) \ dy)^\frac{3b+2}{5} \right ] \\
& \leq & C \left (\Vert \rho(t) \Vert_q R^{-q}(x) \right )^b
\end{eqnarray*}
\noindent for $-\frac{4}{1 + b} < -3$, which is satisfied since $b < \frac{5}{18}$.

Combining the estimates for $I$, $II$, and $III$, the lemma follows.\\

\noindent {\bf \underline{Proof} (Lemma $3$) :} \ Let $q,T > 0$ be
given with $\Vert \rho(t) \Vert_\infty \leq C$ and $\Vert \nabla
\rho(t) \Vert_\infty \leq C$ for all $t \in [0,T]$.  Let $a \in
[0,1)$ be given with $a \leq \frac{3}{q}$. Consider $\vert x \vert \geq 1$ and define
$$\eta := \left (\Vert \rho(t) \Vert_q R^{-q}(x) \right )^a.$$

For any $i,k = 1,2,3$, we have

\begin{equation*}
\begin{array}{lll}
\vert \partial_{x_i} E_k(t,x) \vert & = & \vert \int_{\vert x - y
\vert < \eta} \partial_{y_i} \rho(t,y) \frac{(x-y)_k}{\vert x - y
\vert^3} \ dy \vert + \vert \int_{\vert x - y \vert = \eta}
\rho(t,y) \frac{(x-y)_k}{\vert x - y \vert^3} \frac{(x-y)_i}{\vert
x - y \vert} \ dS_y \vert \\
\\
& \ & + \vert \int_{\vert x - y \vert > \eta} \rho(t,y)
\partial_{y_i} (\frac{(x-y)_k}{\vert x - y
\vert^3}) \ dy \ \vert \\
\\
& =: & I + II + III.
\end{array}
\end{equation*}

Then,

\begin{eqnarray*}
I & \leq & \int_{\vert x - y \vert < \eta} \Vert \nabla_x \rho(t) \Vert_\infty \vert x - y \vert^{-2} \ dy \\
& \leq & C \eta.
\end{eqnarray*}

We estimate $II$ for the large $\vert x \vert$ and small $\vert x
\vert$ cases.  For $\vert x \vert > 2\eta$,

\begin{eqnarray*}
II & \leq & \int_{\vert x - y \vert = \eta} (\Vert \rho(t) \Vert_q
R^{-q}(y))^a \vert \rho(t,y) \vert^{1-a} \vert x - y \vert^{-2} \ dS_y \\
& \leq & \Vert \rho(t) \Vert_q^a \Vert \rho(t)
\Vert_{L^\infty}^{1-a} \eta^{-2} \int_{\vert x - y \vert =
\eta} R^{-aq}(\vert x \vert - \eta) \ dS_y \\
& \leq & C \Vert \rho(t) \Vert_q^a \ \eta^{-2} \
R^{-aq}(\frac{1}{2}\vert x \vert) \eta^2 \\
& \leq & C (\Vert \rho(t) \Vert_q R^{-q}(x))^a.
\end{eqnarray*}

For $\vert x \vert < 2\eta$,

\begin{eqnarray*}
II & \leq & \Vert \rho(t) \Vert_\infty \eta^{-2} \int_{\vert x - y \vert = \eta} dS_y \\
& \leq & C \\
& \leq & C (2\eta \vert x \vert^{-1}) \\
& \leq & C \eta.
\end{eqnarray*}

Then,
\begin{eqnarray*}
III & \leq & C \int_{\eta < \vert x - y \vert < \frac{1}{2} \vert
x \vert} \vert \rho(t,y) \vert \ \vert x - y \vert^{-3} \ dy + C
\int_{\vert x - y \vert > \frac{1}{2} \vert x \vert} \vert
\rho(t,y) \vert \ \vert x - y \vert^{-3} \ dy \\
& =: & A + B.
\end{eqnarray*}

Estimating $A$, we have for $n \in [0,1)$,

\begin{eqnarray*}
A & \leq & \Vert \rho(t) \Vert_q^n \int_{\eta < \vert x - y \vert
< \frac{1}{2} \vert x \vert} \vert \rho(t,y) \vert^{1-n} R^{-nq}
(y) \vert x - y \vert^{-3} \ dy \\
& \leq & C \Vert \rho(t) \Vert_q^n R^{-nq} (\frac{1}{2} \vert x
\vert) \int_{\eta < \vert x - y \vert < \frac{1}{2} \vert x \vert}
(k^{\frac{1 - n}{2}}
(t,y) + k^\frac{3(1-n)}{5}(t,y) ) \vert x - y \vert^{-3} \ dy \\
& \leq & C ( \Vert \rho(t) \Vert_q R^{-q}(x))^n \ \left [ \
(\int_{\eta < \vert x - y \vert < \frac{1}{2} \vert x \vert} \vert
x - y \vert^{-\frac{6}{1+n}} \ dy )^\frac{1+n}{2} \right . \\
& \ & \hspace{1.5in} + \ \left . (\int_{\eta < \vert x - y \vert <
\frac{1}{2} \vert x \vert} \vert x - y
\vert^{-\frac{15}{3n+2}})^\frac{3n+2}{5} \ \right ] \\
& \leq & C ( \Vert \rho(t) \Vert_q R^{-q}(x) )^n \ \left [
(\int_\eta^\infty r^{-\frac{6}{1+n} + 2} \ dr)^\frac{1+n}{2} +
(\int_\eta^\infty r^{-\frac{15}{3n+2} + 2} \ dr)^\frac{3n+2}{5}
\right ] \\
& \leq & C ( \Vert \rho(t) \Vert_q R^{-q}(x) )^n \ (
\eta^{\frac{3}{2}(1+n) - 3} + \eta^{\frac{3}{5}(3n+2) -3}) \\
& \leq & C ( \Vert \rho(t) \Vert_q R^{-q}(x) )^n \
(\eta^{-\frac{3}{2}(1-n)} + \eta^{-\frac{9}{5}(1-n)}) \\
& \leq & C ( \Vert \rho(t) \Vert_q R^{-q}(x) )^n  \ \left \{
\begin{array}{cc}  \eta^{-\frac{9}{5}(1-n)}, & \eta \leq 1 \\
\eta^{-\frac{3}{2}(1-n)}, & \eta \geq 1 \end{array} \right.
\end{eqnarray*}

\noindent For $\eta \leq 1$, we may choose $n = \frac{14a}{9a + 5}$, and for $\eta \geq 1$, we may
choose $n = \frac{5a}{3a + 2}$. Thus,

$$ A \leq C( \Vert \rho(t) \Vert_q R^{-q}(x) )^a. $$

Finally, we estimate $B$ and find
\begin{eqnarray*}
B & \leq & \int_{\vert x - y \vert > \frac{1}{2} \vert x \vert} \vert \rho(t,y) \vert \
\vert x - y \vert^{-3} \ dy \\
& \leq & \Vert \rho(t) \Vert_q^a (\frac{1}{2} \vert x \vert)^{-aq}
\int_{\vert x - y \vert > \frac{1}{2} \vert x \vert} \vert
\rho(t,y) \vert^{1-a} \vert x - y \vert^{aq-3}
R^{-aq}(y) \ dy \\
& \leq & C ( \Vert \rho(t) \Vert_q R^{-q}(x) )^a \int_{\vert x - y \vert > \frac{1}{2} \vert x \vert} (k^\frac{1-a}{2}(t,y) +
k^{\frac{3}{5}(1-a)}(t,y)) \ R^{-3}(y) \ dy \\
& \leq & C ( \Vert \rho(t) \Vert_q R^{-q}(x) )^a \left [ (\int R^{-\frac{6}{1+a}}(y) \ dy)^\frac{1+a}{2} + (\int R^{-\frac{15}{3a+2}}(y)
\ dy)^\frac{3a+2}{5} \right ] \\
& \leq & C ( \Vert \rho(t) \Vert_q R^{-q}(x) )^a
\end{eqnarray*}
\noindent for $-\frac{6}{1 + a} < -3$, which is satisfied since $a < 1$.\\

Combining the estimates for $I$, $II$, $A$, and $B$, the lemma follows.\\

\noindent {\bf \underline{Proof} (Lemma $4$) :} \  Let $T > 0$ and $q \in [0, \frac{54}{13})$ be given with
$\Vert \rho(t) \Vert_\infty \leq C$ and $\Vert \nabla \rho(t) \Vert_\infty \leq C$ for all $t \in [0,T]$
and $Q_f(T) < \infty$.  Let $t \in [0,T]$ be given.
For any $D > 0$, taking $\vert x \vert \leq D$, we find
\begin{equation}
\label{smallx}
\begin{array}{rcl}
\vert \rho(t,x) \vert & \leq & \Vert \rho(t) \Vert_\infty \\
\\
& \leq & \Vert \rho(t) \Vert_\infty \frac{R^q(D)}{R^q(x)} \\
\\
& \leq & C R^{-q}(x).
\end{array}
\end{equation}

Now, define $C^{(3)} := \max \{1, 8TQ_g(T) \}$ and let $\vert x
\vert \geq C^{(3)}$. Define for every $t \in [0,T]$ and $x \in
\mathbb{R}^3$, $$ \mathcal{E}(t,x) = E(t,x) + A(t,x).$$

From the Vlasov equation,
$$\frac{\partial}{\partial s} \left ( g(s,X(s), V(s)) \right ) = - \mathcal{E}(s,X(s)) \cdot \nabla_v F(V(s)), $$
and thus
\begin{equation}
\label{vlasovg}
g(t,x,v) = g(0,X(0),V(0)) - \int_0^t \mathcal{E}(s,X(s)) \cdot \nabla_v F(V(s)) \ ds.
\end{equation}

Thus, to estimate $\rho$, we must consider $\int \mathcal{E}(s,X(s)) \cdot \nabla_v F(V(s)) \ dv$.
Assume $f$ is nonzero along $(X(s), V(s))$.  Then,

\begin{eqnarray*}
\left \vert \int \mathcal{E}(s,X(s)) \cdot \nabla_v F(V(s)) dv
\right \vert & \leq &
\left \vert \int \mathcal{E}(s,X(s)) \cdot (\nabla_v F(V(s)) - \nabla_v F(v + \int_s^t \mathcal{E}(\tau, x) \ d\tau)) \ dv \right \vert \\
& \ & + \ \left \vert \int \left( \mathcal{E}(s,X(s)) -
\mathcal{E}(s, x + (s-t)v) \right) \right. \\
& \ & \left. \ \ \ \ \ \cdot \nabla_v F(v + \int_s^t \mathcal{E}(\tau, x) \ d\tau) \ dv \right \vert \\
& \ & + \ \left \vert \int  \nabla_v \cdot (F(v) \ \mathcal{E}(s,x + (s-t)v)) \ dv \right \vert \\
& \ & + \ \left \vert \int  F(v) \nabla_v \cdot (\mathcal{E}(s,x + (s-t)v)) \ dv \right \vert \\
& =: & I + II + III + IV.
\end{eqnarray*}

Using Lemmas $2$ and $3$, as well as (III), we find $a \in [0,1)$
and $b \in [0,\frac{5}{18})$ with $a \leq \frac{3}{q}$, $b \leq
\frac{2}{q}$, and $ a + b = 1$ such that
$$ \vert \mathcal{E}(t,x) \vert \leq C \left ( \Vert \rho(t) \Vert_q R^{-q}(x) \right )^b$$
and
$$ \vert \nabla \mathcal{E}(t,x) \vert \leq C \left ( \Vert \rho(t) \Vert_q R^{-q}(x) \right )^a.$$

By the Mean Value Theorem, for $\tau \in [s,t]$ and $i = 1,2,3$,
there exist $\xi_1^i$ on the line segment between $X(\tau)$ and
$x$ such that
$$ \mathcal{E}_i(\tau,X(\tau)) - \mathcal{E}_i(\tau,x) = \nabla_x \mathcal{E}_i(\tau, \xi_1^i) \cdot (X(\tau) - x).$$

Hence,
\begin{eqnarray*}
I & \leq & \int_{\vert v \vert \leq Q_g(T)} \vert \mathcal{E}(s,X(s)) \vert \ \Vert \nabla^2 F \Vert_\infty \
\vert V(s) - \left ( v + \int_s^t \mathcal{E}(\tau,x) \ d\tau \right ) \vert \ dv \\
& \leq &  C \int_{\vert v \vert \leq Q_g(T)} \left ( \Vert \rho(s) \Vert_q R^{-q}(X(s)) \right )^b \
\left \vert \int_s^t ( \mathcal{E}(\tau, X(\tau)) - \mathcal{E}(\tau,x) ) \ d\tau \right \vert \ dv \\
& \leq & C \int_{\vert v \vert \leq Q_g(T)} \left ( \Vert \rho(s) \Vert_q R^{-q}(X(s)) \right )^b \
\int_s^t \sup_i \vert \nabla \mathcal{E}(\tau, \xi_1^i) \vert \ \vert X(\tau) - x \vert \ d\tau \ dv \\
& \leq & C \int_{\vert v \vert \leq Q_g(T)} \left ( \Vert \rho(s)
\Vert_q R^{-q}(X(s)) \right )^b \ T Q_g(T) \int_s^t \sup_i \left (
\Vert \rho(\tau) \Vert_q R^{-q}(\xi_1^i) \right )^a \ d\tau \ dv.
\end{eqnarray*}

Since we know, by (\ref{xchar1}), for any $i = 1,2,3$

\begin{eqnarray*}
\vert \xi_1^i \vert & \geq & \vert X(\tau) \vert - \vert X(\tau) - x \vert \\
& \geq & \frac{1}{2} \vert x \vert - T Q_g(T) \\
& \geq & \frac{1}{4} \vert x \vert,
\end{eqnarray*}

we find

\begin{eqnarray*}
I & \leq & C  \left ( \Vert \rho(s) \Vert_q R^{-q}(\frac{1}{2}
\vert x \vert) \right )^b R^{-aq}(\frac{1}{4} \vert x \vert)
\int_s^t \Vert \rho(\tau) \Vert_q^a \ d\tau \\
& \leq & C R^{-q(a+b)}(x) \left ( \Vert \rho(s) \Vert_q^b \int_s^t \Vert \rho(\tau) \Vert_q^a \ d\tau \right ).
\end{eqnarray*}

Similarly, using the above lemmas and the Mean Value Theorem, for
any $i = 1,2,3$ there exist $\xi_2^i$ between $X(\tau)$ and $x +
(s-t)v$ such that
\begin{eqnarray*}
II & \leq & C \sum_{i=1}^3 \int_{\vert v \vert \leq Q_g(T)} \sup_i \vert \nabla
\mathcal{E}(s,\xi_2^i) \vert \ \vert X(s) - (x + (s-t)v) \vert \ \Vert \nabla F \Vert_\infty \ dv \\
& \leq & C \sum_{i=1}^3 \int_{\vert v \vert \leq Q_g(T)} \sup_i \left ( \Vert
\rho(s) \Vert_q R^{-q}(\xi_2^i) \right )^a \ \vert \int_s^t
\int_\tau^t \mathcal{E}(\iota, X(\iota)) \ d\iota \ d\tau \vert \
dv.
\end{eqnarray*}

Since we know, using (\ref{xchar1}), for any $i = 1,2,3$

\begin{eqnarray*}
\vert \xi_2^i \vert & \geq & \vert X(s) \vert - \vert X(s) - (x + (s-t)v) \vert \\
& \geq & \frac{1}{2} \vert x \vert - \vert \int_s^t ( V(\tau) - v ) \ d\tau \vert \\
& \geq & \frac{1}{2} \vert x \vert - 2TQ_g(T) \\
& \geq & \frac{1}{4} \vert x \vert,
\end{eqnarray*}

it follows that

\begin{eqnarray*}
II & \leq & C \int_{\vert v \vert \leq Q_g(T)} \left ( \Vert \rho(s) \Vert_q R^{-q}(\frac{1}{4} \vert x \vert) \right )^a
\ \int_s^t \left ( \Vert \rho(s) \Vert_q R^{-q}(X(\tau)) \right )^b \ \ d\tau \ dv \\
& \leq & C R^{-q(a+b)}(x) \left ( \Vert \rho(s) \Vert_q^a \ \int_s^t \Vert \rho(\tau) \Vert_q^b \ d\tau \right ).
\end{eqnarray*}

By the Divergence Theorem,

$$ III = 0, $$

and finally,

\begin{eqnarray*}
IV & \leq & \int_{\vert v \vert \leq Q_g(T)} \Vert F \Vert_\infty \vert \nabla_x \cdot \mathcal{E}(s, x + (s-t)v) \vert \ \vert s-t \vert \ dv \\
& \leq & C \int_{\vert v \vert \leq Q_g(T)} \vert \rho(s, x + (s-t)v) \vert \ dv \\
\end{eqnarray*}

Collecting the estimates for $I - IV$, we have

\begin{equation}
\label{edotf}
\begin{array}{rcl} \vert \int \mathcal{E}(s,X(s)) \cdot \nabla_v F(V(s)) dv \vert & \leq &
C \left ( R^{-q(a+b)}(x) ( \Vert \rho(s) \Vert_q^b \int_s^t \Vert \rho(\tau) \Vert_q^a \ d\tau) \right. \\
\\
& \ & \ \ + R^{-q(a+b)}(x) ( \Vert \rho(s) \Vert_q^a \ \int_s^t \Vert \rho(\tau) \Vert_q^b \ d\tau) \\
\\
& \ & \ \ \left. + \int_{\vert v \vert \leq Q_g(T)} \vert \rho(s, x + (s-t)v) \vert \ dv \right ). \end{array}
\end{equation}

Since $a + b = 1$, we proceed from (\ref{vlasovg}) and using (\ref{edotf}), $(IV)$, and (\ref{vsupport}), we find
\begin{eqnarray*}
\vert \rho(t,x) \vert & = & \vert \int g(t,x,v) \ dv \vert \\
& \leq & \int_{\vert v \vert \leq Q_g(T)} \vert g(0,X(0),V(0)) \vert \ dv
+ \vert \int_0^t \int_{\vert v \vert \leq Q_g(T)} \mathcal{E}(s,X(s)) \cdot \nabla_v F(V(s)) \ dv \ ds \vert \\
& \leq & CR^{-q}(x) + C \int_0^t \left ( R^{-q(a+b)}(x)
( \Vert \rho(s) \Vert_q^b \int_s^t \Vert \rho(\tau) \Vert_q^a \ d\tau) \right. \\
& \ & \ \ \left. + R^{-q(a+b)}(x) ( \Vert \rho(s) \Vert_q^a \ \int_s^t \Vert \rho(\tau) \Vert_q^b \ d\tau) +
\int_{\vert v \vert \leq Q_g(T)} \vert \rho(s, x + (s-t)v) \vert \ dv \right ) \ ds \\
& \leq & C R^{-q}(x) \left ( 1 + \int_0^t \Vert \rho(s) \Vert_q^b \int_s^t \Vert \rho(\tau) \Vert_q^a \ d\tau ds +
\int_0^t \Vert \rho(s) \Vert_q^a \ \int_s^t \Vert \rho(\tau) \Vert_q^b \ d\tau ds  \right ) \\
& \ & \ \ + C \int_{\vert v \vert \leq Q_g(T)} \int_0^t \vert \rho(s,x + (s-t)v) \vert \ ds \ dv \\
& \leq & C R^{-q}(x) \left ( 1 + (\int_0^t \Vert \rho(s) \Vert_q^b ds ) \ (\int_0^t \Vert \rho(\tau) \Vert_q^a \ d\tau) +
(\int_0^t \Vert \rho(s) \Vert_q^a \ ds) \ (\int_0^t \Vert \rho(\tau) \Vert_q^b \ d\tau)  \right ) \\
& \ & \ \ + C \int_{\vert v \vert \leq Q_g(T)} \int_0^t \vert \rho(s,x + (s-t)v) \vert \ ds \ dv.
\end{eqnarray*}
Then, applying H\"{o}lder's inequality twice,
\begin{eqnarray*}
\vert \rho(t,x) \vert & \leq & C R^{-q}(x) \left ( 1 + 2T^2 \left
( \int_0^t \Vert \rho(s) \Vert_q \ ds \right )^{a + b} \right. \\
& \ & \ \ \ \left. + R^q(x) \int_{\vert v \vert \leq Q_g(T)} \int_0^t \Vert \rho(s) \Vert_q R^{-q}(x + (s-t)v) \ ds \ dv \right ) \\
& \leq & C R^{-q}(x) \left ( 1 + \int_0^t \Vert \rho(s) \Vert_q \ ds +
R^{q}(x) \ R^{-q}(\frac{7}{8} \vert x \vert) \ Q^3_g(T) \ \int_0^t \Vert \rho(s)\Vert_q \ ds \right ) \\
& \leq & C R^{-q}(x) \left (1 + \int_0^t \Vert \rho(s) \Vert_q \ ds \right ).
\end{eqnarray*}
Combining this with (\ref{smallx}) (with $D = C^{(3)}$) and multiplying by $R^q(x)$, we have for all $x$,
$$ R^q(x) \vert \rho(t,x) \vert \leq C (1 + \int_0^t \Vert \rho(s) \Vert_q \ ds ).$$
Since the right side of the inequality is independent of $x$, we take the supremum over all $x$ to find,
$$ \Vert \rho(t) \Vert_q \leq C (1 + \int_0^t \Vert \rho(s) \Vert_q \ ds) $$
and by Gronwall's Inequality,
$$ \Vert \rho(t) \Vert_q \leq C. $$

Thus, $\Vert \rho(t) \Vert_q$ is bounded for all $t \in [0,T]$.  Notice, too, that the choice of $a + b = 1$ forces
$q < \frac{54}{13}$,
and the proof of the lemma is complete.\\

\appendix

\section{Appendix A}

In this appendix, we will explore the argument used in Section $4.2.6$ of \cite{Glassey} to bound derivatives of $f$ and $E$.

Let $D$ be any $x$ derivative.  Then, using the Vlasov equation,
$$ \partial_t(Df) + v \cdot \nabla_x(Df) - E \cdot \nabla_v(Df) =
DE \cdot \nabla_v f$$ so that $$ \frac{d}{ds} Df(s, X(s), V(s)) =
DE \cdot \nabla_v f(s, X(s), V(s)).$$ Hence, $$ \vert Df(t,x,v)
\vert \leq \vert Df(0, X(0), V(0)) \vert + \int_0^t \vert DE \cdot
\nabla_v f(s, X(s), V(s)) \vert \ ds.$$ Define
$$ \vert f(s) \vert_1 = \sup_{x,v} \vert \partial_x f(s,x,v) \vert
+ \sup_{x,v} \vert \partial_v f(t,x,v) \vert$$ and
$$ \vert E(s) \vert_1 = \sup_x \vert \partial_x E(s,x) \vert.$$
Then, $$ \vert Df(t,x,v) \vert \leq C + \int_0^t \vert E(s)
\vert_1 \vert f(s) \vert_1 \ ds.$$ We see that $\partial_v f$
satisfies an inequality of similar type because
$$ \partial_t (\partial_v f) + v \cdot \nabla_x (\partial_v f) - E
\cdot \nabla_v (\partial_v f) = - \partial_v v \cdot \nabla_x f.$$
It follows that $$ \vert f(t) \vert_1 \leq C + \int_0^t ( 1 +
\vert E(s) \vert_1) \vert f(s) \vert_1 \ ds.$$ However, from
(\ref{derivfield}), we can conclude (with $d = \Vert \nabla_x \rho \Vert_\infty)$ that $$ \vert E(s) \vert_1 \leq C (1 + \ln^* \vert f(s) \vert_1 )$$ where
$$ \ln^*s = \left \{ \begin{array}{lr} s & 0 \leq s \leq 1 \\
1 + \ln s & s > 1. \end{array} \right.$$ Therefore, for $t \leq T$,
$$ \vert f(t) \vert_1 \leq C \left ( 1 + \int_0^t (1 + \ln^* \vert
f(s) \vert_1 ) \vert f(s) \vert_1 \ ds. \right ).$$ It now follows
by an application of Gronwall's Iequality that $$ \vert f(t) \vert_1
\leq C$$ and $$ \vert E(t) \vert_1 \leq C$$ for $t \leq T$. This
concludes the argument to bound field derivatives and derivatives of
the density, and thus ends Appendix A.

\section{Appendix B}
The bound on $Q(t)$ is obtained as follows.  From (\ref{char}) and (\ref{GBU}), we have
\begin{eqnarray*}
\vert \hat{V}(t) \vert & \leq & \vert \hat{V}(t - \Delta) \vert +
\int_{t - \Delta}^t \vert E(s,\hat{X}(s)) \vert \ ds \\
& \leq & Q(t - \Delta) + C\Delta Q^\frac{13}{15}(t).
\end{eqnarray*}

Since the above constant is independent of the particular
characteristic, we find $$ Q(t) \leq Q(t - \Delta) + C \Delta
Q^\frac{13}{15}(t)$$ for
\begin{eqnarray*}
\Delta & = & \min \left \{t, \frac{1}{4C^{(2)}}
Q^{-\frac{4}{3}}(t) \cdot Q^\frac{13}{20}(t) \right \}. \\
& = & \min \left \{t, \frac{1}{4C^{(2)}} Q^{-\frac{41}{60}}(t)
\right \}
\end{eqnarray*}
Since $Q$ is non-decreasing, there exists $T_1$ such that
$$ \Delta = \left \{ \begin{array}{cc} t & t \leq T_1 \\
\frac{1}{4C^{(2)}} Q^{-\frac{41}{60}} (t)& t \geq T_1. \end{array} \right.$$
Take $t_0$ in the interval of existence.  Without loss of generality,
$t_0 > T_1$. Let $$ t_1 = t_0 - \frac{1}{4C^{(2)}}
Q^{-\frac{41}{60}}(t_0),$$
$$ t_{i+1} = t_i - \frac{1}{4C^{(2)}} Q^{-\frac{41}{60}}(t_i) \ \ \ (i =
1,2,...)$$ as long as $t_i > T_1$. Then
$$ t_i - t_{i+1} = \frac{1}{4C^{(2)}} Q^{-\frac{41}{60}}(t_i) \geq \frac{1}{4C^{(2)}}
Q^{-\frac{41}{60}}(t_0)$$ which is a uniform lower bound on the
length of each subinterval.  So, there is a first $i$, say $i =
k$, such that $t_k \leq T_1$.  Thus, $t_k \geq 0$ and therefore
$$ (t_0 - t_1) + (t_1 - t_2) + \cdots + (t_{k-1} - t_k) \geq k
\cdot \frac{1}{4C^{(2)}} Q^{-\frac{41}{60}}(t_0)$$ which implies
that
$$ Q^{-\frac{41}{60}}(t_0) \cdot k \leq 4C^{(2)}t_0.$$
Now we have
\begin{eqnarray*}
Q(t_0) & = & Q(t_k) + \sum_{i=0}^{k-1} [Q(t_i) - Q(t_{i+1})] \\
& \leq & Q(t_k) + C \sum_{i=0}^{k-1} \Delta \cdot Q^\frac{13}{15}(t_i) \\
& \leq & Q(T_1) + C  \sum_{i=0}^{k-1} \frac{1}{4C^{(2)}}
Q^{-\frac{41}{60}}(t_0) \cdot Q^\frac{13}{15}(t_i) \\
& \leq & Q(T_1) + C \cdot k Q^{-\frac{41}{60}}(t_0) \cdot
Q^\frac{13}{15} (t_0) \\
& \leq & C t_0 Q^\frac{13}{15}(t_0).
\end{eqnarray*}
Therefore, $Q(t_0)$ is bounded, and the proof is complete.  This
ends Appendix B.


\begin{thebibliography}{99}
\bibitem{B} Batt, J. Global symmetric solutions of the
initial-value problem of stellar dynamics. J. Diff. Eq. {\bf 1977}, 25 :
342-364.

\bibitem{BR} Batt, J.; Rein, G. Global classical solutions of
the periodic Vlasov-Poisson system in three dimensions. C. R.
Academy of Sci. {\bf 1991}, 313(1): 411-416.

\bibitem{C} Caglioti, E.; Caprino, S.; Marchioro, C.; Pulvirenti, M.
The Vlasov equation with infinite charge. Arch. Ration. Mech.
Anal. {\bf 2001}, 159:85-108.

\bibitem{Glassey} Glassey, R. {\it The Cauchy Problem in Kinetic
Theory}; S.I.A.M: Philadelphia, 1996.

\bibitem{GS} Glassey, R.; Strauss, W. Singularity formation in a
collisionless plasma could occur only at high velocities.
Arch. Ration. Mech. Ana. {\bf 1986}, 92:59-90.

\bibitem{Horst} Horst, E. On the classical solutions of the
initial value problem for the unmodified nonlinear
Vlasov-equation, part I. Math. Methods Appl. Sci. {\bf 1981}, 3:229-248.

\bibitem{Horst82} Horst, E. On the classical solutions of the
initial value problem for the unmodified nonlinear
Vlasov-equation, part II. Math. Methods Appl. Sci. {\bf 1982},
4:19-32.

\bibitem{Horst93} Horst, E. On the asymptotic growth of the
solutions of the Vlasov-Poisson system. Math. Methods Appl. Sci.
{\bf 1993}, 16:75-85.

\bibitem{J} Jabin, P. E. The Vlasov-Poisson system with
infinite charge and energy. J. Statist. Phys. {\bf 2001},
103(5/6):1107-1123.

\bibitem{K} Kurth, R. Das anfangswertproblem der
stellardynamik. Z. Astrophys. {\bf 1952} 30:213-229.

\bibitem{LP} Lions, P.L.; Pertham, B. Propogation of moments
and regularity for the three dimensional Vlasov-Poisson system.
Invent. Math. {\bf 1991}, 105:415-430.

\bibitem{OU} Okabe, S.; Ukai, T. On classical solutions in the
large in time of two-dimensional Vlasov's equation. Osaka J. Math.
{\bf 1978}, 15:245-261.

\bibitem{Pankavich} Pankavich, S. Global existence for the
Radial Vlasov-Possion System with Steady Spatial Asymptotics.
{\bf 2004} \ (submitted for publication)

\bibitem{P} Pfaffelmoser, K. Global classical solution of the
Vlassov-Poisson system in three dimensions for general initial
data. J. Diff. Eq.  {\bf 1992}, 95(2):281-303.

\bibitem{RR} Rein, G.; Rendall, A. Global existence of classical
solutions to the Vlasov-Poisson system in a three dimensional
cosmological setting. Arch. Ration. Mech. Anal. {\bf 1994}
126:183-201.

\bibitem{Schaeffer} Schaeffer, J. Global existence of smooth solutions
to the Vlasov-Poisson system in three dimensions. Comm. PDE. {\bf
1991}, 16(8/9):1313-1335.

\bibitem{SSAVP} Schaeffer, J. Steady spatial asymptotics
for the Vlasov-Poisson system. Math. Methods Appl. Sci. {\bf
2003}, 26:273-296.

\bibitem{VPSSA} Schaeffer, J. The Vlasov-Poisson system
with steady spatial asymptotics. Comm. PDE. {\bf 2003},
28(5/6):1057-1084.

\bibitem{W} Wollman, S. Global-in-time solutions of the two
dimensional Vlasov-Poisson system. Comm. Pure Appl. Math. {\bf
1980}, 33:173-197.


\end{thebibliography}
\end{document}